\documentclass[12pt]{article}
\usepackage{amssymb, amsthm, amsfonts, amsmath, amscd}
\usepackage{enumitem}
\usepackage[utf8]{inputenc}
\usepackage[margin=1in]{geometry}
\usepackage{graphicx}
\usepackage{tikz}
\usepackage{tikz-cd}
\usepackage{mathrsfs}
\usepackage{hyperref}
\usepackage{array}
\usepackage{multirow}
\usepackage{ytableau}
\usepackage{graphicx}
\usepackage{lineno}
\usepackage{color}
\usepackage{comment}

\newtheorem{thm}{Theorem}

\newtheorem{proposition}[thm]{Proposition}

\theoremstyle{definition}

\newtheorem{definition}[thm]{Definition}

\newtheorem{example}[thm]{Example}
\newtheorem*{problem1}{Problem 1}
\newtheorem*{problem2}{Problem 2}
\newtheorem*{problem3}{Problem 3}
\newtheorem*{problem2a}{Problem 2$'$}
\newtheorem*{problem3a}{Problem 3$'$}

\usepackage{lipsum}
\usepackage{tikz}
\usepackage{comment}
\usepackage{tikz-cd}
\usepackage{mathrsfs}

\newcommand{\Y}{\mathbb Y}

\newcommand{\R}{\mathbb R}

\newcommand{\Z}{\mathbb Z}

\newcommand{\F}{\mathbb F}
\newcommand{\De}{\Delta}
\newcommand{\la}{\lambda}
\newcommand{\vp}{\varphi}
\newcommand{\om}{\omega}
\newcommand{\Sym}{\operatorname{Sym}}
\newcommand{\HL}{\operatorname{HL}}
\newcommand{\Mat}{\operatorname{Mat}}
\newcommand{\Harm}{\operatorname{Harm}}

\newcommand{\nil}{0}
\newcommand{\id}{\operatorname{id}}
\newcommand{\wt}{\widetilde}
\newcommand{\GL}{\textrm{GL}}
\newcommand{\UU}{\textrm U}
\newcommand{\CC}{\mathcal C}
\newcommand{\calF}{\mathcal F}
\renewcommand{\L}{\mathcal L}
\newcommand{\U}{\mathcal U}
\renewcommand{\P}{\mathcal P}
\newcommand{\G}{\mathcal G}
\newcommand{\g}{\mathfrak{g}}
\newcommand{\gl}{\mathfrak{gl}}
\renewcommand{\u}{\mathfrak u}

\newcommand{\T}{\mathcal T}
\newcommand{\Lat}{\mathbb L}
\newcommand{\GLB}{\mathbb{GLB}}
\newcommand{\UB}{\mathbb{UB}}
\newcommand{\E}{\mathscr{E}}
\newcommand{\Up}{\uparrow}
\newcommand{\Down}{\downarrow}
\newcommand{\naB}{\wt{\nabla}}
\newcommand{\DeB}{\wt{\Delta}}
\newcommand{\bfla}{\pmb{\la}}
\newcommand{\bfmu}{\pmb{\mu}}

\newcommand{\M}{\mathcal M}

\title{Mackey-type identity for invariant functions on Lie algebras of finite unitary groups and an application}

\author{Cesar Cuenca$^a$ and Grigori Olshanski$^{b,c,d}$\\
\\
{\small $^a$Department of Mathematics, Harvard University, Cambridge, MA, USA}\\
{\small Email: cesar.a.cuenk@gmail.com}\\
{\small $^b$Institute for Information Transmission Problems, Moscow, Russia}\\
{\small $^c$Skolkovo Institute of Science and Technology, Moscow, Russia}\\
{\small $^d$Faculty of Mathematics, HSE University, Moscow, Russia}\\
{\small Email: olsh2007@gmail.com}}
\date{}

\begin{document}

\maketitle

\centerline{Dedicated to Karl H. Hofmann on the occasion of his 90th birthday}

\abstract{
The Mackey-type identity mentioned in the title relates the operations of parabolic induction and restriction for invariant functions on the Lie algebras of the finite unitary groups $U(N, \F_{q^2})$. This result is applied to constructing positive harmonic functions on a new branching graph with a negative Hall-Littlewood parameter, as introduced in the authors' paper [Adv. Math. vol. 395 (2022), 108087].
This in turn implies the existence of an infinite-parameter family of invariant measures for the coadjoint action of an infinite-dimensional analogue of the groups $U(N, \F_{q^2})$.
}

\section{Preface}

In the authors' previous paper \cite{CO}, a novel branching graph was introduced and studied. The vertices of this graph are the Young diagrams with an even number of boxes\footnote{A second branching graph was defined in \cite{CO}, where the vertices are Young diagrams of odd size, but we restrict ourselves to the even branching graph in this note.} and the edges are endowed with certain positive weights (``formal multiplicities'') deriving from the Hall-Littlewood symmetric functions with a negative parameter. An open problem proposed in that paper is the classification of \emph{positive harmonic functions} on this branching graph. In the present paper, we construct an infinite-parameter family of positive harmonic functions. In the process we produce a Lie-algebra analogue of the classical Mackey's theorem (\cite[Theorem 1]{Mackey}, \cite[\S44]{CR}) which is of independent interest.

The open problem above is equivalent to the classification of Radon measures, invariant with respect to the coadjoint action of an infinite-dimensional unitary group over a finite field $\F_{q^2}$. In this guise, our proposed problem is a new addition to the vast literature of asymptotic representation theory. For example, the analogous question for the infinite-dimensional unitary group $U(\infty)$ (the field is $\R$ and not a finite field) was studied by Olshanski-Vershik \cite{OV}. Moreover, a close relative is the classification problem of characters and spherical functions for infinite-dimensional matrix groups.
Work in this direction was done by Thoma \cite{T} and Vershik-Kerov \cite{VK1} for the infinite symmetric group; by Voiculescu \cite{V}, Boyer \cite{B}, Vershik-Kerov \cite{VK2},  and Okounkov-Olshanski \cite{OO1}, \cite{OO2} for infinite-dimensional classical Lie groups; by Gorin-Olshanski \cite{GOl} in the quantized setting; by Gorin-Kerov-Vershik \cite{GKV} for the infinite-dimensional general linear group over $\F_q$; and many others.

Note that the invariant Radon measures from our paper \cite{CO} are infinite measures. In this respect our theory drastically differs from Fulman’s probabilistic theory of random matrices over finite fields (see his survey paper \cite{F} and references therein). However, in both cases, Hall–Littlewood symmetric functions play a fundamental role.

In terms of techniques, our first idea is to view invariant measures (equivalently, harmonic functions on branching graphs) as linear functionals on the space of invariant functions on Lie algebras of unitary groups. The space of invariant functions has the structure of a \emph{twisted} bimodule; a similar ``twist'' was observed before by van Leeuwen \cite{vL}, who studied representations of finite groups of Lie type. By contrast, the space of invariant functions on general linear Lie algebras is a relative of the representation ring studied by Zelevinsky \cite{Z}, and it is a bialgebra.
The proofs of the bimodule/bialgebra structures depend on a version of Mackey's theorem for parabolic induction and restriction in the context of Lie algebras. For this part, we drew ideas from the works of van Leeuwen and of Zelevinsky.

Another important idea in our work is the modification of a construction, due to Sergei Kerov, of functionals on the algebra of symmetric functions that are positive on Macdonald functions, see \cite[Sec. 4]{GO} for an exposition. Kerov's construction depends substantially on the coproduct of the algebra of symmetric functions; we manage to adapt his idea by using instead the coaction map of our twisted bimodule.

The present note is an announcement; detailed proofs will be given in another paper.

\section{Invariant Radon measures}\label{sec_2}

Consider one of the classical series $\{G(n)\}$ of finite groups of Lie type. They form a nested chain under natural inclusions, thus giving rise to the inductive limit $G(\infty) := \varinjlim{G(n)}$. We further take certain topological completion $\overline{G} \supset G(\infty)$ to obtain an infinite-dimensional topological group $\overline{G}$ which is locally compact and separable.

The example with most history comes from the general linear groups $G(n) = \GL(n, \F_q)$ over a finite field $\F_q$ of cardinality $q$. This classical series gives rise to $G(\infty) = \GL(\infty, \F_q)$. The completion $\overline{G}$ is denoted by $\GLB$ and consists of infinite $\Z_{\ge 1}\times \Z_{\ge 1}$ matrices that are almost upper-triangular (have finitely many nonzero entries below the diagonal) and are invertible. The group $\GLB$ was discovered and studied by Vershik and Kerov starting in 1982, see \cite{V82}, \cite{VK98}, \cite{VK07}; their ideas were completed and further expanded by Gorin-Kerov-Vershik \cite{GKV}. It was argued by these authors that $\GLB$ is the ``right'' infinite-dimensional analogue of the $\GL(n, \F_q)$'s because of its rich representation theory (unlike the case of $\GL(\infty, \F_q)$) and its connections to the infinite symmetric group.

The completion $\overline{G}$ can also be constructed for the classical series of unitary, orthogonal and symplectic groups over a finite field and a natural question, motivated by Gorin-Kerov-Vershik \cite{GKV}, is to classify the semifinite traces (characters) of these infinite-dimensional classical groups. In this note, we are interested in a Lie algebra analogue of this problem. Namely, instead of  characters or representations of $\overline{G}$, we are dealing with \emph{$\overline{G}$-invariant Radon measures} for the coadjoint action. Note that for $\overline{G}$, there exist natural notions of the (infinite-dimensional) Lie algebra $\overline{\mathfrak{g}}$ and the dual space  $\overline{\mathfrak{g}}^*$, as well as of adjoint and coadjoint actions on $\overline{\mathfrak{g}}$ and $\overline{\mathfrak{g}}^*$, respectively. Furthermore, $\overline{\mathfrak{g}}^*$ is a locally compact vector space with respect to a natural topology, which is consistent with the coadjoint action of $\overline{G}$. So the notion of invariant Radon measures on $\overline{\mathfrak{g}}^*$ makes sense.\footnote{Recall that an infinite measure on a locally compact space is Radon if it is finite on all compact subsets.} We remark that the $\overline{G}$-invariant Radon measures on $\overline{\mathfrak{g}}^*$ turn out to be infinite measures, with the obvious exception of the delta measure at the zero matrix.

\begin{problem1}
Study the Radon measures on $\overline{\mathfrak{g}}^*$ that are invariant with respect to the coadjoint action of $\overline{G}$. In particular, classify the ergodic invariant Radon measures.
\end{problem1}

When $\overline{G} = \mathbb{GLB}$, the dual space $\overline{\mathfrak{g}}^*$ is denoted by $\Lat(q)$; it consists of all $\Z_{\ge 1}\times\Z_{\ge 1}$ matrices with finitely many nonzero entries below and on the diagonal, and the coadjoint action of $\GLB$ is given by conjugation. Problem 1 was completely solved by the authors in this case, see \cite[Sec. 4]{CO}, though the actual difficult part of the argument is due to Matveev \cite{Mat}. However, Problem 1 is still open for groups $\overline{G}$ coming from the other classical series.

In this note, we concentrate on the unitary group setting. Whenever we discuss finite unitary groups, we switch from the field $\F_q$ to its quadratic extension $\F_{q^2}$; we also assume that $q$ is odd. For any matrix $A = (a_{i, k})$ with entries in $\F_{q^2}$, define
$$
A^{[q]} :=(a_{i, k}^q),\qquad A^* := (A^{[q]})^T = (A^T)^{[q]},
$$
where the superscript $T$ means transposition. For any $m\in\Z_{\ge 1}$, the finite unitary group $\UU(m, \F_{q^2})$ is by definition the group of invertible matrices $g\in\GL(m, \F_{q^2})$ such that
\begin{equation}\label{unitary_condition}
g^* J_m\, g = J_m,
\end{equation}
where $J_m$ is the following $m\times m$ matrix:
\begin{equation}\label{matrix_J}
J_m := \begin{bmatrix}  0 & &  &  & 1\\ & \ddots & & \ddots & \\ & &1& & \\ & \ddots & & \ddots & \\ 1 & & & & 0 \end{bmatrix}.
\end{equation}

It turns out that in order to construct the completion $\overline{G}$, we must examine separately finite unitary groups of even or odd dimension, thus leading to two parallel theories. For simplicity, we restrict ourselves to the even dimension case here, so that $G(n) = \UU(2n, \F_{q^2})$.
The inclusions $\UU(2n, \F_{q^2}) \hookrightarrow \UU(2n+2, \F_{q^2})$ are the ones obtained naturally by labeling the rows and columns of each matrix in $\UU(2m, \F_{q^2})$ by $-m, \cdots, -1, 1, \cdots, m$; they lead to the inductive limit $G(\infty) = \UU(2\infty, \F_{q^2})$.
The completion $\overline{G}$ is denoted here by $\mathbb{UB}$; it consists of infinite matrices with rows and columns parametrized by $\Z\setminus\{0\}$, with finitely many nonzero entries below the diagonal, and satisfying the infinite analogue of the unitarity condition \eqref{unitary_condition}.
(In our previous paper \cite{CO}, we used the superscripts $\E$ and $\mathscr O$ to distinguish between the ``even'' and ``odd'' cases; in particular, the group $\mathbb{UB}$ was denoted there by $\mathbb{UB}^\E$.)

The Lie algebra $\overline{\mathfrak{g}}$ of $\mathbb{UB}$ consists of $(\Z\setminus\{0\})\times (\Z\setminus\{0\})$ matrices $X$ with finitely many entries below the diagonal, and satisfying the skew-Hermitian condition 
$$
X^*J_{2\infty} + J_{2\infty}X = 0,
$$
where $J_{2\infty}$ is the $(\Z\setminus\{0\})\times (\Z\setminus\{0\})$ matrix with $1$'s on the second diagonal, as in \eqref{matrix_J}. The dual space $\overline{\mathfrak{g}}^*$, denoted by $\wt\Lat(q^2)$, is defined similarly but with the addition that matrices have finitely many nonzero entries also on the diagonal. The latter condition allows one to define a pairing $\overline{\mathfrak{g}}\times \overline{\mathfrak{g}}^*\to \F_q$,
making $\overline{\mathfrak{g}}^*$ the dual space to $\overline{\mathfrak{g}}$.

In the remainder of this note, our aim is to study Problem 1 in the unitary case, when $\overline{G} = \UB$ and $\overline{\mathfrak{g}}^* = \wt{\Lat}(q^2)$. Even though the general linear case ($\overline{G} = \GLB$ and $\overline{\mathfrak{g}}^* = \Lat(q)$) is well-understood by now, we also discuss it for comparison and because it will help us obtain results for the unitary setting.

\section{Branching graphs and harmonic functions}

In this section, we translate Problem 1 to the language of harmonic functions on branching graphs.
We also explain our results from \cite{CO} that connect the branching graphs to the Hall-Littlewood symmetric functions.

\subsection{Finitely additive functions on cylinder sets}

In all instances of the general setting, there is a natural notion of \emph{cylinder sets} in $\overline{\mathfrak{g}}^*$. These sets are open and compact, and form a ring of sets as well as a base of the topology of $\overline{\mathfrak{g}}^*$. Any Radon measure induces a finitely-additive measure on the ring of cylinder sets.
Conversely, any such finitely-additive measure can be extended, in a unique way,  to a Radon measure. The ring of cylinder sets of $\overline{\mathfrak{g}}^*$ is invariant with respect to the coadjoint action of $\overline{G}$. Thus,  the convex cone of invariant Radon measures on $\overline{\mathfrak{g}}^*$ is affine-isomorphic to the convex cone of invariant, finitely-additive nonnegative functions on the ring of cylinder sets of $\overline{\mathfrak{g}}^*$.

\smallskip

In the general linear group case $\overline{G} = \GLB$, cylinder sets are defined as follows. An \emph{elementary cylinder set of level $n\in\Z_{\ge0}$} in $\overline{\mathfrak{g}}^* = \Lat(q)$ is composed of the matrices $M = [m_{i, j}]_{i,,j=1}^\infty$ which have a prescribed $n\times n$ top-left corner $[m_{i, j}]_{i, j = 1}^n$ and satisfy the triangularity condition $m_{i,j}=0$, for all $(i,j)$ such that $i>n$ and $i\ge j$.
An arbitrary cylinder set of level $n$ is a disjoint union (necessarily finite) of elementary cylinder sets of the same level. Any cylinder set of level $n$ is also a cylinder set of level $n+1$. In particular, each elementary cylinder set of level $n$ decomposes into a disjoint union of some elementary cylinder sets of level $n+1$. This fact is the base of the relation \eqref{gl_harm} below. Finally, we say that two elementary cylinder sets of level $n$ are \emph{equivalent} if their $n\times n$ top-left corners  are conjugate by a matrix from $\GL(n, \F_q)$. Thus, the equivalence classes of level $n$ are parametrized by the set of $\GL(n, \F_{q})$-orbits in $\gl(n, \F_q)$; we denote this set by $\T_n$ and form the disjoint union
\begin{equation}\label{GL_orbits}
\T := \bigsqcup_{n\ge 0}{\T_n}.
\end{equation}

If $F : \T\to\R_{\ge 0}$ is the restriction of a $\GLB$-invariant measure on $\Lat(q)$ to the (equivalence classes of) elementary cylinder sets, the constraint of finite-additivity translates to the following relations:
\begin{equation}\label{gl_harm}
F(\bfmu) = \sum_{\bfla\in\T_{n+1}}{ L^{n+1}_n(\bfla, \bfmu) F(\bfla) },\qquad \bfmu\in\T_n,\quad n=0, 1, 2, \cdots.
\end{equation}
The quantities $L^{n+1}_n(\bfla, \bfmu)$, by definition, are
\begin{equation}\label{L_link}
L^{n+1}_n(\bfla, \bfmu) := \#\left\{ x\in(\F_q)^n \ \ \Bigg| \ \begin{bmatrix} X_{\bfmu} & x \\ 0 & 0 \end{bmatrix} \text{ belongs to the $\GL(n\!+\!1,\, \F_q)$-orbit $\bfla$} \right\},
\end{equation}
where $X_{\bfmu}$ is an arbitrary matrix belonging to the orbit $\bfmu$. The conclusion is that the convex cone of $\GLB$-invariant measures from Problem 1 is isomorphic to the convex cone of functions $F: \T\to\R_{\ge 0}$ satisfying the harmonicity condition \eqref{gl_harm}. The ergodic invariant measures correspond to the \emph{extreme} positive harmonic functions, i.e. to those in the extreme rays of the convex cone.

\smallskip

In the unitary group case $\overline{G} = \UB$, cylinder sets of $\overline{\mathfrak{g}}^* = \wt\Lat(q^2)$ are defined in a similar way. A difference now is that matrices grow in both directions. According to this, $n\times n$ corners are replaced by central submatrices $[m_{i, j}]_{i, j = -n, \cdots, -1, 1, \cdots, n}$ of order $2n$. Such a submatrix necessarily belongs to $\u(2n, \F_{q^2})$. The equivalence classes of elementary cylinder sets are parametrized by
$$\wt\T := \bigsqcup_{n\ge 0}{\wt\T_{2n}},$$
where $\wt\T_{2n}$ is the set of $\UU(2n, \F_{q^2})$-orbits in $\u(2n, \F_{q^2})$.

If $H : \wt\T\to\R_{\ge 0}$ is the restriction of a $\UB$-invariant measure on $\wt\Lat(q^2)$ to the elementary cylinder sets, the finite-additivity property is equivalent to the relations
\begin{equation}\label{u_harm}
H(\wt\bfmu) = \sum_{\wt\bfla\in\wt\T_{2n+2}}{ K^{n+1}_n(\wt\bfla, \wt\bfmu)\, H(\wt\bfla) },\qquad \wt\bfmu\in\wt\T_{2n},\quad n = 0, 1, 2, \cdots,
\end{equation}
where
\begin{equation*}
\begin{aligned}
K^{n+1}_n(\wt\bfla, \wt\bfmu) := \#\left\{ \!(x, y)\in (\F_{q^2})^{2n} \times \F_{q^2} \ \Bigg| \right.&
\,y^q = -y, \text{ and } \begin{bmatrix} 0 & -x^*J_{2n} & y\\ 0 & X_{\wt\bfmu} & x \\ 0 & 0 & 0 \end{bmatrix} \text{ belongs} \\
&\qquad\quad \left.\text{ to the $\UU(2n\!+\!2,\, \F_{q^2})$-orbit $\wt\bfla$} \right\};
\end{aligned}
\end{equation*}
here, $X_{\wt\bfmu}$ is any matrix in $\u(2n, \F_{q^2})$ belonging to $\wt\bfmu$. The conclusion is the same: $\UB$-invariant measures are in a bijective correspondence with nonnegative functions on the set of unitary orbits $\wt\T$, satisfying the relations  \eqref{u_harm}.

\subsection{Branching graphs}

The relations \eqref{gl_harm} and \eqref{u_harm} above can be interpreted as the harmonicity conditions on certain branching graphs. We next explain this terminology.

By a \emph{branching graph} we mean a pair $\Gamma=(V,\tau)$, where $V$ is the vertex set of a graph and $\tau$ is a strictly positive function on its edges. The graph is graded, meaning that its vertices are partitioned into countably many levels $V = V_0 \sqcup V_1 \sqcup V_2 \sqcup \cdots$, and such that the endpoints of any edge belong to adjacent levels. We further assume:

$\bullet$ Each level $V_n$ is finite, and $V_0$ consists of a single element called the \emph{root vertex}.

$\bullet$ For each $v_n\in V_n$, there exists at least one $v_{n+1}\in V_{n+1}$ such that $\{v_{n+1}, v_n\}$ is an edge. Likewise, if $n\ne 0$, there exists at least one $v_{n-1}\in V_{n-1}$ such that $\{v_n, v_{n-1}\}$ is an edge.

We regard $\tau(v_{n+1}, v_n)$, the value of $\tau$ on the edge $\{v_{n+1}, v_n\}$, as its \emph{weight} or else as its \emph{formal multiplicity}. Let us emphasize that this quantity is not necessarily an integer.

It is convenient to set $\tau(v_{n+1}, v_n) = 0$, whenever $\{v_{n+1}, v_n\}$ is not an edge; then the set of edges can be deduced from the values of $\tau$.

\begin{example}\label{young}
Let $\Y_n$ be the set of partitions of size $n$ and $\Y = \Y_0 \sqcup \Y_1 \sqcup \Y_2 \sqcup \cdots$ be the set of all partitions. As in Macdonald \cite{Mac}, we identify partitions with their Young diagrams. The size of a diagram $\la\in\Y$ is denoted by $|\la|$. If $\mu, \la\in\Y$ are such that $\mu\subset\la$ and $\la\setminus\mu$ is a single box, then we write $\mu\nearrow\la$; observe that this implies $|\la| - |\mu| = 1$. 
The \emph{Young graph}, see e.g. \cite{BO}, is the branching graph with graded vertex set $\Y$ and the edges $\mu\nearrow\la$. All edge multiplicities are equal to $1$: $\tau(\la, \mu) = 1$ whenever $\mu\nearrow\la$.
\end{example}

For a branching graph $\Gamma = (V, \tau)$, we say that a real-valued function $F : V\to\R$ is a \emph{positive harmonic function on $\Gamma$} if:

$\bullet$ (Positivity) $F(v) \ge 0$, for all $v\in V$.

$\bullet$ (Harmonicity) For any $n=0,1,2,\cdots$, and $v\in V_n$, we have
\begin{equation}\label{general_harm}
F(v) = \sum_{w\in V_{n+1}}{\tau(w, v) F(w)}.
\end{equation}
The set of positive harmonic functions on $\Gamma$ is a convex cone that we denote $\Harm_{\ge 0}(\Gamma)$.
By definition, the \emph{(minimal) boundary of $\Gamma$} is the union of extreme rays of the convex cone $\Harm_{\ge 0}(\Gamma)$.

Note that any branching graph $\Gamma=(V,\tau)$ is connected, because each vertex $v\in V$ is joined by a path to the root vertex. We also need a slight extension of the above definitions to the case of disconnected graphs. By a \emph{disconnected branching graph} we mean a disjoint union of a countable number of branching graphs. We also allow shifts of gradings on the connected components, so that the root vertex of a component can now be assigned another level --- not necessarily $0$, but some positive integer. The definition of harmonic functions, whose main component is the relation \eqref{general_harm}, remains the same for disconnected branching graphs. The set of harmonic functions is still a convex cone, and in fact it is isomorphic to the product of the cones corresponding to the connected components. The boundary is a disjoint union of the boundaries of the connected components.

\begin{definition}
The \emph{graph of\/ $\GLB$-invariant measures} $\Gamma^{\GLB}$ is the graph with vertices $\T = \T_0\sqcup\T_1\sqcup\T_2\sqcup\cdots$, such that an edge joins $\bfmu\in\T_n$ and $\bfla\in\T_{n+1}$ if and only if $L^{n+1}_n(\bfla, \bfmu) \ne 0$; the weight of that edge is set to $\tau(\bfla, \bfmu) := L^{n+1}_n(\bfla, \bfmu)$. One can show that $\Gamma^{\GLB}$ is a disconnected branching graph.
\end{definition}

\begin{definition}
The \emph{graph of $\UB$-invariant measures} $\Gamma^{\UB}$ is the graph with vertices $\wt\T = \wt\T_0\sqcup\wt\T_2\sqcup\wt\T_4\sqcup\cdots$, and such that an edge joins $\wt\bfmu\in\wt\T_{2n}$ and $\wt\bfla\in\wt\T_{2n+2}$ if and only if $K^{n+1}_n(\wt\bfla, \wt\bfmu)\ne 0$. The edge-weights are defined to be $\tau(\wt\bfla, \wt\bfmu) := K^{n+1}_n(\wt\bfla, \wt\bfmu)$. By definition, the vertices of $\wt\T_{2n}$ are placed on level $n$. Again, one can show that $\Gamma^{\UB}$ is a disconnected branching graph.
\end{definition}

Let, as above, $\overline{G}$ be one of the groups $\GLB$, $\UB$. We denote by $\M^{\overline{G}}$ the cone of $\overline{G}$-invariant Radon measures on $\overline{\mathfrak{g}}^*$. The general harmonicity relation \eqref{general_harm} specializes to \eqref{gl_harm} and \eqref{u_harm} for the graphs $\Gamma^{\overline{G}}=\Gamma^{\GLB}$ and $\Gamma^{\overline{G}}=\Gamma^{\UB}$, respectively.  Thus, we have the isomorphism $\M^{\overline{G}} \cong \Harm_{\ge 0}(\Gamma^{\overline{G}})$, and the ergodic $\overline{G}$-invariant measures correspond to the extreme rays of the cone $\Harm_{\ge 0}(\Gamma^{\overline{G}})$, i.e. to the boundary of $\Gamma^{\overline{G}}$. As a result, Problem 1 is equivalent to the following.

\begin{problem2} 
Study the convex cones of positive harmonic functions $\Harm_{\ge 0}(\Gamma^{\GLB})$ and $\Harm_{\ge 0}(\Gamma^{\UB})$. Find the boundary of the disconnected branching graphs $\Gamma^{\GLB}$ and $\Gamma^{\UB}$.
\end{problem2}

\subsection{The nilpotent part and Hall-Littlewood symmetric functions}\label{sec_nilpotent}

In the general setting, there is an analogue of nilpotency for infinite matrices in $\overline{\g}^*$: \emph{pro-nilpotent matrices}. The set of $\overline{G}$-invariant Radon measures supported on the pro-nilpotent matrices of $\overline{\g}^*$ will be denoted $\M_0^{\overline{G}}$. This is a subcone of $\M^{\overline{G}}$, corresponding to a distinguished connected component $\Gamma_0^{\overline{G}}$ of the disconnected branching graph $\Gamma^{\overline{G}}$.

It is known \cite{CO} that the cone $\M^{\GLB}$ is isomorphic to the product of countably many copies of $\M_0^{\GLB}$, and therefore determining the structure of $\M^{\GLB}_0$ would not be a simplification of Problem 1, but rather an equivalent question.

For the group $\UB$ the picture is similar but with an interesting modification. Recall that $\UB$ is our shorthand notation for the ``even'' version $\UB^\E$; there is also an ``odd'' version $\UB^{\mathscr O}$. It turns out that the cone $\M^{\UB^{\E}}$ is isomorphic to a doubly infinite product of cones, with countably many copies of $\M^{\UB^\E}_0$ and countably many copies of $\M^{\UB^{\mathscr O}}_0$. Thus, Problem 1 for the group $\UB$ reduces to the study of the cones $\M^{\UB^\E}_0$ and  $\M^{\UB^{\mathscr O}}_0$. These two cases are similar; below we focus on the case of $\M^{\UB^\E}_0$.

\subsubsection{The general linear case: $\Gamma^{\GLB}_0\sim \Y^{\HL}(q^{-1})$}\label{sec_qdeformed}

Let us describe the vertices and edge-multiplicities of the branching graph $\Gamma^{\overline{G}}_0$ when $\overline{G} = \GLB$.
The vertices of $\Gamma^{\GLB}_0$ are divided into levels, where the $n$-th level consists of the nilpotent orbits within $\T_n$; these are parametrized by partitions $\la\in\Y_n$.
In fact, an $n\times n$ nilpotent matrix belongs to the orbit $\la$ if and only if its Jordan normal form is composed of blocks of sizes $\la_1, \la_2, \cdots$. Next, if $\la, \mu$ are partitions of sizes $n+1,\, n$, respectively, then $\{\la, \mu\}$ is an edge of $\Gamma^{\GLB}_0$ if and only if the quantity $L^{n+1}_n(\la, \mu)$ defined in \eqref{L_link} is nonzero. It was proved by Kirillov \cite[Sec. 2.3]{K} and Borodin \cite[Thm. 2.3]{B}, see also \cite[Sec. 9]{CO}, that $L^{n+1}_n(\la, \mu) \ne 0$ if and only if $\mu\nearrow\la$. Moreover, if $\mu\nearrow\la$ and the single box $\la\setminus\mu$ is in column $k$, then
\begin{equation}\label{L_cotransition}
L^{n+1}_n(\la, \mu) = \begin{cases}
q^{n - \sum_{j\ge k}{m_j(\mu)}}(1 - q^{-m_{k-1}(\mu)}) & \text{if }k>1,\\
q^{n - \sum_{j\ge 1}{m_j(\mu)}} & \text{if }k=1,
\end{cases}
\end{equation}
where  $m_j(\mu): = \#\{ i \mid \mu_i = j\}$, for any $j\in\Z_{\ge 1}$.

An important fact is that the branching graph $\Gamma^{\GLB}_0$ just described has a connection with the Hall-Littlewood symmetric functions. In order to make this connection, let us recall some of the theory of symmetric functions from Macdonald \cite{Mac}.

Denote by $\Sym$ the real algebra of symmetric functions, which is generated by $1$ and the algebraically independent power sums $p_1, p_2, \cdots$. We will need the \emph{Hall-Littlewood symmetric functions} (HL functions for short), which are parametrized by partitions, form a basis of $\Sym$, and depend on a parameter $t$, which for us will belong to $(0, 1)$. There are two versions of the HL functions, which differ by a normalization constant, namely the $P$ and $Q$ versions. We work exclusively with the $Q$-HL functions which are denoted by $Q_{\nu}(; t)$, $\nu\in\Y$; see \cite[Ch. III]{Mac} for their precise definition.

The algebra $\Sym$ is graded by declaring $\deg(p_k) = k$, for all $k$, and $\deg(1)=0$. The HL function $Q_{\nu}(; t)$ is homogeneous of degree $|\nu|$, for any $\nu\in\Y$. Then for any $n=0, 1, 2, \cdots$ and $\mu\in\Y_n$, we have an expansion
$$p_1\cdot Q_{\mu}(; t) = \sum_{\la\in\Y_{n+1}}{\psi_{\la/\mu}(t)\, Q_{\la}(; t)},$$
with  some coefficients $\psi_{\la/\mu}(t)$. An explicit formula exists: $\psi_{\la/\mu}(t) \ne 0$ if and only if $\mu\nearrow\la$ and in this case, if $k$ is the column number of $\la\setminus\mu$, then
\begin{equation}\label{psi_cotransition}
\psi_{\la/\mu}(t) = \begin{cases}
1-t^{m_{k-1}(\mu)} & \text{if }k>1,\\
1 & \text{if }k=1.\end{cases}
\end{equation}
Observe that $\psi_{\la/\mu}(t)\ge 0$, for any $t\in (0, 1)$, and in particular for $t = q^{-1}$. Thus we are allowed to make the following definition.

\begin{definition}
The \emph{HL-deformed Young graph} $\Y^{\HL}(q^{-1})$ is the branching graph with the same vertices and edges as the usual Young graph from Example \ref{young}, but the edges $\mu\nearrow\la$ have formal multiplicities $\kappa(\la, \mu) := \psi_{\la/\mu}(q^{-1})$.
\end{definition}

Observe that the branching graphs $\Gamma^{\GLB}_0$ and $\Y^{\HL}(q^{-1})$ have the same vertices and edges. Moreover by comparing \eqref{L_cotransition} and \eqref{psi_cotransition}, the edge-multiplicities are related by
$$
L^{n+1}_n(\la, \mu) = \psi_{\la/\mu}(q^{-1})\cdot
\frac{q^{\sum_{i\ge 1}{(i-1)\mu_i} - {|\mu| \choose 2}}}{q^{\sum_{i\ge 1}{(i-1)\la_i} - {|\la| \choose 2}}},
\qquad \mu\in\Y_n,\quad \la\in\Y_{n+1},\quad \mu\nearrow\la.
$$
This means that the branching graphs $\Gamma^{\GLB}_0$ and $\Y^{\HL}(q^{-1})$ are \emph{similar} in the sense of Kerov \cite[Sec. 4]{Ke0}; we denote this by $\Gamma^{\GLB}_0\sim \Y^{\HL}(q^{-1})$. As a result, the convex cones $\Harm_{\ge 0}(\Gamma^{\GLB}_0)$ and $\Harm_{\ge 0}(\Y^{\HL}(q^{-1}))$ are isomorphic; one can go from one cone to the other by multiplying by the function $f(\la) = q^{\sum_{i\ge 1}{(i-1)\la_i} - {|\la| \choose 2}},\ \la\in\Y$.

\subsubsection{The unitary case: $\Gamma^{\UB}_0\sim \Y^{\HL}_{\E}(-q^{-1})$}

Let us explain a similar connection between the branching graph $\Gamma^{\UB}_0$ and HL functions with a negative parameter.

The vertices of $\Gamma^{\UB}_0$ are partitioned into levels, where the $n$-th level is the set of nilpotent orbits of $\wt\T_{2n}$.
As it turns out, the $\UU(2n, \F_{q^2})$-orbit of a nilpotent matrix $N\in\u(2n, \F_{q^2})$ is determined by its type, i.e. by the $\GL(2n, \F_{q^2})$-orbit of $N$ when regarded as a matrix in $\gl(2n, \F_{q^2})$. Hence the $n$-th level of the vertex-set of $\Gamma^{\UB}_0$ is the set of partitions of size $2n$.

As for the edge-multiplicities $K^{n+1}_n(\la, \mu)$, $\la\in\Y_{2n+2},\, \mu\in\Y_{2n}$, of the branching graph $\Gamma^{\UB}_0$, we need a new notation. If two partitions $\mu, \la$ are such that $\mu\subset\la$, and $\la\setminus\mu$ consists exactly of two boxes on the same column, or two boxes in adjacent columns, then we write $\mu\!\nearrow\!\!\nearrow\!\la$. Then $K^{n+1}_n(\la, \mu) \ne 0$ if and only if $\mu\!\nearrow\!\!\nearrow\!\la$; moreover in this latter case, an explicit formula for $K^{n+1}_n(\la, \mu)$ was obtained in \cite[Prop. 8.2]{CO}.

\smallskip

In order to relate $K^{n+1}_n(\la, \mu)$ to the HL functions, let us replace the HL parameter $t\in(0, 1)$ by $-t\in(-1, 0)$ because we want to make it clear that we are working with a negative parameter. Consider the following normalization of $Q$-HL functions:
$$\wt{Q}_{\nu}(; -t)\!:=\!(-1)^{\sum_{i\ge 1}{\!(i-1)\nu_i}}\!\cdot Q_{\nu}(; -t),\quad \nu\in\Y.$$
The set $\{\wt{Q}_{\nu}: \nu\in\Y\}$ is a basis of $\Sym$, and each $\wt{Q}_{\nu}$ is homogeneous of degree $|\nu|$.
Consequently, for any partition $\mu\in\Y_n$, $n=0,1,2,\cdots$, we have an expansion of the form
\begin{equation}\label{branch_two}
((1 - t^2)p_2)\cdot\wt{Q}_{\mu}(; -t) = \!\sum_{\la\in\Y_{n+2}}{\xi_{\la/\mu}(-t)\,\wt{Q}_{\la}(; -t)},
\end{equation}
for some coefficients $\xi_{\la/\mu}(-t)$. One can verify that $\xi_{\la/\mu}(-t)\ge 0$, for any $-t\in (-1, 0)$, and $\xi_{\la/\mu}(-t)\ne 0$ if and only if $\mu\!\nearrow\!\!\nearrow\!\la$; this is true in particular for $-t = -q^{-1}$.

\begin{definition}
The \emph{even HL-deformed Young graph}, denoted $\Y^{\HL}_{\mathscr{E}}(-q^{-1})$, has graded set of vertices $\Y_{\mathscr{E}} := \Y_0\sqcup\Y_2\sqcup\Y_4\sqcup\cdots$, and an edge joins $\mu\in\Y_{2n}$ and $\la\in\Y_{2n+2}$ if and only if $\mu\!\nearrow\!\!\nearrow\!\la$. The edge $\mu\!\nearrow\!\!\nearrow\!\la$ has multiplicity $\kappa(\mu, \la) := \xi_{\la/\mu}(-q^{-1})$.
\end{definition}

One of the main results from \cite{CO} is Theorem 8.9 there, which shows in particular that $\Gamma^{\UB}_0$ and $\Y^{\HL}_{\mathscr{E}}(-q^{-1})$ are similar, meaning that the edge-multiplicities $K^{n+1}_n(\la, \mu)$ and $\xi_{\la/\mu}(-q^{-1})$ differ by a multiplicative gauge factor (of the form $f(\la)/f(\mu)$). Consequently, $\Harm_{\ge 0}(\Gamma^{\UB}_0) \cong \Harm_{\ge 0}(\Y^{\HL}_{\mathscr{E}}(-q^{-1}))$.

\section{Interlude}

Let us summarize our discussion so far.

\medskip

$\bullet$ The proposed Problem 1 is to study the convex cone $\M^{\overline{G}}$ of $\overline{G}$-invariant Radon measures on the topological vector space $\overline{\g}^*$, and in particular the ergodic measures. This problem is the origin of this note and of our previous paper \cite{CO}. 

\smallskip

$\bullet$ By inspecting the countable topological base of cylinder sets of $\overline{\g}^*$, one obtains an isomorphism between $\M^{\overline{G}}$ and the cone of positive harmonic functions on certain disconnected branching graph $\Gamma^{\overline{G}}$. With this terminology, our question is equivalent to Problem 2. The set of ergodic measures within $\M^{\overline{G}}$ corresponds to the boundary of $\Gamma^{\overline{G}}$.

\smallskip

$\bullet$ The measures supported on pro-nilpotent matrices form a distinguished convex subcone $\M^{\overline{G}}_0 \subset\M^{\overline{G}}$, and the understanding of its structure would essentially solve Problem 1 (and 2). This question also has a translation to the language of branching graphs and asks for the boundary of a distinguished connected component $\Gamma^{\overline{G}}_0$ of $\Gamma^{\overline{G}}$.

\smallskip

$\bullet$ In the cases of our interest, when $\overline{G}$ is $\GLB$ or $\UB$, it was discovered in \cite{CO} that the branching graphs $\Gamma^{\GLB}_0$ and $\Gamma^{\UB}_0$ are similar to $\Y^{\HL}(q^{-1})$ and $\Y^{\HL}_{\mathscr{E}}(-q^{-1})$, respectively; these last two are built from HL functions with parameters $q^{-1}$ and $-q^{-1}$, respectively. Hence Problems 1 and 2 get simplified to the following:

\begin{problem3}
 Study the convex cones $\Harm_{\ge 0}(\Y^{\HL}(q^{-1}))$ and $\Harm_{\ge 0}(\Y^{\HL}_{\mathscr{E}}(-q^{-1}))$. Find the boundary of the branching graphs $\Y^{\HL}(q^{-1})$ and $\Y^{\HL}_{\mathscr{E}}(-q^{-1})$.
\end{problem3}

More generally, we can ask about $\Harm_{\ge 0}(\Y^{\HL}(t))$ and $\Harm_{\ge 0}(\Y^{\HL}_{\mathscr{E}}(-t))$, for any $t\in(0, 1)$, but we can give them a Lie theoretic interpretation only when $t = q^{-1}$ and $q$ is the power of a prime number (an odd prime number, in the unitary case).

\smallskip

$\bullet$ Both $\Y^{\HL}(q^{-1})$ and $\Y^{\HL}_{\mathscr{E}}(-q^{-1})$ resemble other well-known branching graphs of representation theoretic origin that are defined from the Pieri rule of certain bases of the ring of symmetric functions, see e.g. \cite{BO}, \cite{KOO}, and references therein. This is significant because there are many techniques in the literature to study such branching graphs.
In fact, this connection to HL functions was essential to solve Problem 1 (and 2) in the case of general linear groups by drawing on the knowledge of $\Y^{\HL}(q^{-1})$, as provided by a special case of a result of Matveev \cite{Mat}. However, the study of the unitary version $\Y^{\HL}_{\mathscr{E}}(-q^{-1})$ and its boundary seem to require new ideas to cope with the negative HL parameter and with the power sum $p_2$ in \eqref{branch_two}, which are new features.

\medskip

We now give a preview of the remaining of this note. Our goal is to address Problem 2 (equivalently, Problem 1) for the unitary case. Ultimately, we would like to describe the boundaries of $\Gamma^{\UB}$ and $\Gamma^{\UB}_0$, and more generally to understand the structure of $\Harm_{\ge 0}(\Gamma^{\UB})$ and $\Harm_{\ge 0}(\Gamma^{\UB}_0) \cong \Harm_{\ge 0}(\Y^{\HL}_{\mathscr{E}}(-q^{-1}))$. We do not solve this problem completely, but obtain a partial result, namely Theorem \ref{kerov_construction}, which shows how to produce an infinite-parameter family of examples of positive harmonic functions in $\Harm_{\ge 0}(\Gamma^{\UB})$ (or $\Harm_{\ge 0}(\Gamma^{\UB}_0)$) from just a single one.

The procedure in Theorem \ref{kerov_construction} is an adaptation of \emph{Kerov's mixing construction}, and is inspired by an idea of Sergei Kerov from 1992, who used it in the context of the Young branching graph with a deformed multiplicity function depending on the two Macdonald parameters $(q,t)$. (As was later shown by Matveev \cite{Mat}, the mixing construction makes it possible to obtain all harmonic functions for this graph.)

For the mixing construction, we need yet another description of $\M^{\overline{G}}\cong\Harm_{\ge 0}(\Gamma^{\overline{G}})$ as the convex cone of certain functionals $A^{\overline{G}}\to\R$, on some space $A^{\overline{G}}$.
Such description is known in other similar problems of asymptotic representation theory. The best known example is the equivalence between characters of the infinite symmetric group $S(\infty) = \varinjlim{S(n)}$ and linear functionals $\Sym\to\R$ which are nonnegative on Schur functions, see \cite{BO}.
Here $\Sym$ manifests as the representation ring of the symmetric groups, and the Schur functions represent the characters of the finite symmetric groups $S(n)$.
In this setting, the mixing construction hinges on the bialgebra structure of $\Sym$. More specifically, the key property is that the coproduct $\Delta\!: \Sym \to \Sym\otimes\Sym$ is an algebra homomorphism, and this boils down to a particular case of the classical Mackey's theorem relating the operations of induction and restriction for a symmetric group and its two-block Young subgroups.

In our situation, we work with invariant functions on Lie algebras instead of finite group characters, and we need not the ordinary operations of induction and restriction, but parabolic analogues thereof. Then a suitable version of Mackey's formula is given by Theorem \ref{Mackey_identity}. It underlies our version of the mixing construction.

When $\overline{G} = \GLB$, the corresponding space $A^{\overline{G}}=A^{\GLB}$ (denoted by $A_q$) is a Lie algebra analogue of the representation ring of the finite general linear groups studied by Zelevinsky \cite{Z}. We prove in Proposition \ref{thm_structureA} that $A_q$ is a graded bialgebra.

In the unitary case $\overline{G} = \UB$, an interesting effect arises. In Theorems \ref{thm_twisted_0} and \ref{thm_twisted}, we show that the corresponding space $A^{\overline{G}}=A^{\UB}$ (denoted by $B_{q^2}$) is a \emph{twisted $A_{q^2}$-bimodule}, meaning it is both a module and comodule over $A_{q^2}$, and these structures interact in a nonobvious way. The fact that we have a bimodule and not a bialgebra boils down to the following fact: any maximal Levi subalgebra of $\mathfrak{gl}(n, \F_q)$  is the direct sum of two smaller general linear Lie algebras, while a maximal Levi subalgebra of $\mathfrak{u}(2n, \F_{q^2})$ is a direct sum of the form $\mathfrak{gl}(n_1,\F_{q^2})\oplus\mathfrak{u}(n_2,\F_{q^2})$, where $2n_1+n_2=2n$.  A similar effect is present at the level of group representations, see van Leeuwen \cite{vL}.

\smallskip

In the remainder of this text, we discuss the Mackey-type identity in Section \ref{sec_5}; the spaces $A_q,\, B_{q^2}$ of invariant functions on $\gl(n, \F_q)$'s and $\u(n, \F_{q^2})$'s in Sections \ref{sec_6} and \ref{sec_7}, respectively; in the last two sections, we discuss the mixing construction that produces families of positive harmonic functionals $B_{q^2}\to\R$.

\section{The Mackey-type identity}\label{sec_5}

\subsection{Parabolic induction and parabolic restriction}

Let $G$ be a finite group of Lie type and $\G$ be its Lie algebra. The group $G$ acts on its Lie algebra by the adjoint action $\text{Ad}\!: G\to\text{End}(\G)$.

\begin{definition}
An \emph{invariant function on $\G$} is a real valued function that is constant on the orbits of the adjoint action. Denote the vector space of invariant functions on $\G$ by $\CC(\G)$.
\end{definition}

The operations of \emph{parabolic induction} and \emph{parabolic restriction} depend on the choice of a parabolic subgroup $P < G$ with Levi decomposition
$$P = LU.$$
Here, $L$ is the Levi subgroup and $U$ is the unipotent radical.
Let $\P,\,\L,\,\U$ be the Lie algebras of $P,\, L,\, U$, respectively; they are Lie subalgebras of $\G$ and $\P = \L\oplus\U$ as vector spaces. We use the equality $|U|=|\U|$;  for classical groups, it follows from the existence of the Cayley transform.

\begin{definition}\label{par_res}
The \emph{parabolic restriction} $\Down^{\G}_{\L}: \CC(\G)\to\CC(\L)$ is the linear map
$$(\Down^\G_\L\Psi)(Y) := \frac{1}{|\U|} \sum_{Z\in\U}{\Psi(Y + Z)}, \quad \Psi\in\CC(\G),\ Y\in\L.$$
\end{definition}

\begin{definition}\label{par_ind}
The \emph{parabolic induction} $\Up^\G_\L: \CC(\L)\to\CC(\G)$ is the linear map that takes any $\Phi\in\CC(\L)$ as input and obtains an output by the following procedure:

\smallskip

$\bullet$ Extend $\Phi$ to a function $\Phi': \P\to\R$ by composing with the projection $\pi^\P_\L: \P\twoheadrightarrow \L$ from $\P=\L\oplus\U$ onto the first summand.

$\bullet$ Further extend $\Phi'$ to the whole $\G$ by declaring $\Phi'(W) := 0$, for all $W\in\G\setminus\P$.

$\bullet$ Finally, set
$$
(\Up_\L^\G \Phi)(X):=\sum_{g\in [G/P]}\Phi'(\text{Ad}(g^{-1})X), \quad X\in\G,
$$
where $[G/P]$ is an arbitrary set of representatives for the left cosets of $G$ modulo $P$. 
\end{definition}

Our definitions of parabolic restriction and parabolic induction coincide with the definitions of the truncation map $\tau_{\L\subset\P}^{\G}$ in Lehrer~\cite[Def. (3.1), (i)]{Le}, and the Harish-Chandra induction $\rho^{\G}_{\L\subset\P}$ in \cite[Def. (3.1), (ii)]{Le}, respectively. The condition $|U| = |\U|$ ensures that the parabolic induction and restriction are adjoint to each other with respect to natural inner products, see \cite[Lemma 3.2]{Le}. Moreover, Definitions \ref{par_res} and \ref{par_ind} can be considered as Lie algebra analogues of the maps from van Leeuwen~\cite{vL}.

Beginning in the next section, we specialize $G$ to be either the general linear group $\GL(n, \F_q)$ or the even unitary group $\UU(2n, \F_{q^2})$. In these cases, both $G$ and its Lie algebra $\G$ consist of matrices, the adjoint action is given by conjugation $\text{Ad}(g)(X) = gXg^{-1}$, and the Lie bracket on $\G$ is the commutator $[X, Y] = XY - YX$.

\subsection{The setting for the Mackey-type identity}

Instead of trying to state all the assumptions under which the Mackey-type identity holds, let us restrict ourselves to the following two settings where it will be applied.

\smallskip

\textbf{Setting I.} $\, G = \GL(n, \F_q)$ is the general linear group and its Lie algebra $\G = \gl(n, \F_q)$ consists of all $n\times n$ matrices with entries in $\F_q$.

We will consider parabolic induction and restriction with respect to parabolic subgroups of the following form.
Let $i, j\in\Z_{\ge 0}$ be such that $i+j=n$, and regard $n\times n$ matrices as $2\times 2$ block matrices, where rows and columns are partitioned with respect to $n = i+j$. Then define $P_{i, j}$ as the maximal parabolic subgroup with Levi decomposition $P_{i, j} = L_{i, j}\,U_{i, j}$ given by:
\begin{equation}\label{LU_general}
\begin{aligned}
L_{i, j} &:= \left\{ \begin{bmatrix} A & 0 \\ 0 & B \end{bmatrix} \bigg\rvert\ A\in\GL(i, \F_q),\ B\in\GL(j, \F_q) \right\}\!,\\
U_{i, j} &:= \left\{ \begin{bmatrix} 1 & C \\ 0 & 1 \end{bmatrix} \bigg\rvert\ C\in\Mat_{i\times j}(\F_q) \right\}\!.
\end{aligned}
\end{equation}

Since $L_{i, j}$ is isomorphic to $\GL(i, \F_q)\times\GL(j, \F_q)$, the Lie algebra $\L_{i, j}$ of $L_{i, j}$ is isomorphic to $\gl(i, \F_q)\oplus\gl(j, \F_q)$.

\medskip

\textbf{Setting II.} $\, G = \UU(2n, \F_{q^2})$ is the even unitary group. Its Lie algebra $\G = \u(2n, \F_{q^2})$ is the space of all matrices $X\in\gl(2n, \F_{q^2})$ such that $X^*J_{2n} + J_{2n}X = 0$, where $J_{2n}$ is the matrix in \eqref{matrix_J}.

Fix any $i, j\in\Z_{\ge 0}$ such that $i+j=n$, and regard $(2n)\times (2n)$ matrices as $3\times 3$ block matrices, where rows and columns are partitioned by $2n = i+(2j)+i$. Below we shall consider the maximal parabolic subgroup with Levi decomposition $\wt{P}_{i, j} = \wt{L}_{i, j}\wt{U}_{i, j}$, where
\begin{equation}\label{LU_unitary}
\begin{aligned}
\wt{L}_{i, j} &:= \left\{ \begin{bmatrix} A & 0 & 0 \\ 0 & B & 0 \\ 0 & 0 & A' \end{bmatrix} \bigg\rvert\
A\in\GL(i, \F_{q^2}),\ B\in \UU(2j, \F_{q^2}),\ A' = J_i(A^*)^{-1}J_i \right\}\!,\\
\wt{U}_{i, j} &:= \left\{ \begin{bmatrix} I_i & C & D \\ 0 & I_{2j} & C' \\ 0 & 0 & I_i \end{bmatrix} \bigg\rvert\
C' = -J_{2j}C^*J_i,\quad J_iD + D^*J_i = -J_iCJ_{2j}C^*J_i \right\}\!.
\end{aligned}
\end{equation}

The Lie algebra $\wt\L_{i, j}$ of $\wt{L}_{i, j}$ is isomorphic to $\gl(i, \F_{q^2})\oplus\u(2j, \F_{q^2})$ in this case.

\subsection{The Mackey-type identity}

Let $G$ be a finite group of Lie type from one of the two settings just described, and let $\G$ be its Lie algebra.
Let $P = LU$, $P' = L'U'$ be the Levi decompositions of two maximal parabolic subgroups of $G$ that are both of the form \eqref{LU_general} for general linear groups, or both of the form \eqref{LU_unitary} for even unitary groups.
In particular, $L$ and $L'$ contain the same maximal torus of diagonal matrices. Let $W$ be the corresponding Weyl group embedded as a subgroup (of permutations) of $G$.

\smallskip

For each $w\in W$, consider the maximal parabolic subgroups $P_w < L$ and $P_w' < L'$, given by the Levi decompositions
$$P_w = L_wU_w,\qquad P_w' = L_w'U_w',$$
where
\begin{gather*}
L_w := L\cap w^{-1}L'w,\qquad U_w := L\cap w^{-1}U'w,\\
L_w' := L'\cap wLw^{-1},\qquad U_w' := L'\cap wUw^{-1}.
\end{gather*}
Denote the Lie algebras of $L$ and $L'$ by $\L$ and $\L'$, respectively. The Lie algebras of $L_w$ and $L_w'$ will be
$$\L_w = \L\cap w^{-1}\L'w, \qquad \L'_w = \L'\cap w\L w^{-1}.$$
One can then define the following parabolic induction and restriction maps:
\begin{gather*}
\Up^{\L}_{\L_w}: \CC(\L_w)\to\CC(\L),\qquad\Down^{\L}_{\L_w}: \CC(\L)\to\CC(\L_w),\\
\Up^{\L'}_{\L_w'}: \CC(\L_w')\to\CC(\L'),\qquad \Down^{\L'}_{\L_w'}: \CC(\L')\to\CC(\L_w').
\end{gather*}
Finally, let $T_w: \CC(\L_w)\stackrel{\cong}{\rightarrow}\CC(\L'_w)$ be the linear isomorphism
\begin{equation*}
(T_w f)(X) := f(w^{-1}Xw),\quad f\in\CC(\L_w),\ X\in\L_w'.
\end{equation*}

\begin{thm}[Mackey's theorem for invariant functions on Lie algebras]\label{Mackey_identity}
The following identity between linear maps $\CC(\L)\to\CC(\L')$ holds
\begin{equation}\label{mackey1}
\Down^\G_{\L'}\circ\Up^\G_\L=\sum_{w\in [(W\cap L') \backslash W / (W\cap L)]}\Up^{\L'}_{\L'_w}\circ T_w\circ \Down^\L_{\L_w},
\end{equation}
where $w$ ranges over a set of representatives of $(W\cap L') \backslash W / (W\cap L)$.
\end{thm}

There is another Mackey-type identity for invariant functions on Lie algebras, see Lehrer \cite[(2.6)]{Le} and Springer \cite{Sp}, but that one involves different notions of induction and restriction. Our result should be regarded as the Lie algebra analogue of the Mackey's theorem for finite groups of Lie type, see Digne-Michel \cite[Thm. 5.1]{DM} and van Leeuwen \cite[Thm. 2.3.1]{vL}.

\section{The bialgebra $A_q$}\label{sec_6}

The real graded vector space
$$
A_q := \bigoplus_{n=0}^{\infty}{\CC(\gl(n, \F_q))}
$$
can be equipped with the structure of a graded bialgebra by means of the parabolic induction and parabolic restriction, as will be shown next. The bialgebra $A_q$ will be called the \emph{space of invariant functions of the $\gl(n, \F_q)$'s} and it is the Lie algebra analogue of the representation ring of general linear groups, see Zelevinsky \cite{Z}.

Fix any positive integer $n$ and set
$$G_n := \GL(n, \F_q),\qquad \G_n := \gl(n, \F_q).$$
Consider the maximal parabolic subgroup with Levi decomposition $P_{i, j} = L_{i, j}U_{i, j}$ from \eqref{LU_general}.
The Lie algebra $\L_{i, j}$ of $L_{i, j}$ is isomorphic to $\gl(i, \F_q)\oplus\gl(j, \F_q)$, and so $\CC(\L_{i, j})\cong \CC(\gl(i, \F_q))\otimes\CC(\gl(j, \F_q))$. Therefore the parabolic restriction $\Down^{\G_n}_{\L_{i, j}}$ and parabolic induction $\Up^{\G_n}_{\L_{i, j}}$ are maps of the form:
\begin{gather*}
\Down^{\G_n}_{\L_{i, j}} : \CC(\gl(n, \F_q)) \longrightarrow \CC(\gl(i, \F_q))\otimes\CC(\gl(j, \F_q)),\\
\Up^{\G_n}_{\L_{i, j}} : \CC(\gl(i, \F_q))\otimes\CC(\gl(j, \F_q)) \longrightarrow \CC(\gl(n, \F_q)).
\end{gather*}
These operations naturally induce graded linear maps
$$
\Delta: A_q\to A_q\otimes A_q,\qquad\nabla: A_q\otimes A_q\to A_q.
$$
Indeed, $\Delta$ is defined by the fact that its $n$-th graded part
$$
\left.\Delta\right|_{(A_q)_n}: \CC(\gl(n, \F_q)) \to
\bigoplus_{i, j\in\Z_{\ge 0} :\, i + j = n}{\CC(\gl(i, \F_q))\otimes\CC(\gl(j, \F_q))}
$$
is the sum of $\Down^{\G_n}_{\L_{i, j}}$ over pairs $(i, j)\in (\Z_{\ge 0})^2$ such that $i+j = n$.
Likewise, the map $\nabla$ is defined as the graded map whose $n$-th graded part is the sum of $\Up^{\G_n}_{\L_{i, j}}$ over all pairs $(i, j)\in(\Z_{\ge 0})^2$ such that $i+j = n$.

By definition, the $0$-th graded part $(A_q)_0 = \CC(\gl(0, \F_q))$ is a $1$-dimensional vector space which we identify with $\R$. Let us denote by
$$
e: \R\to A_q,\qquad e^*: A_q\to\R,
$$
the inclusion and projection maps.

\begin{proposition}[The bialgebra structure of $A_q$]\label{thm_structureA}
The maps $\nabla, \De, e, e^*$ make $A_q$ into a bialgebra, which is graded, commutative and cocommutative.
\end{proposition}

This is essentially a known statement, see e.g. \cite[Sec. 2]{Sc} where the bialgebra structure of the subspace of nilpotently supported functions is explained. Also, \cite[Sec. 4]{GR} discusses the bialgebra structure on the space of class functions of general linear groups and its connection to the Hall algebra.

An important part of Proposition \ref{thm_structureA} is the compatibility of the maps $\nabla, \De$, i.e. the fact that $\De$ is an algebra morphism. The proof of this fact closely follows the approach of Zelevinsky \cite{Z} and depends on the Mackey-type identity in Theorem \ref{Mackey_identity}.

\section{The twisted $A_{q^2}$-bimodule $B_{q^2}$}\label{sec_7}

By employing our general paradigm with parabolic induction and restriction, it will turn out that the real graded vector space
$$B_{q^2} := \bigoplus_{n=0}^{\infty}{\,\CC(\u(2n, \F_{q^2}))}$$
is not a graded bialgebra, but rather a \emph{graded twisted $A_{q^2}$-bimodule}; see Definition \ref{def_twist} below. The grading of $B_{q^2}$ is such that $\CC(2n, \F_{q^2})$ is its $n$-th graded part (and not its $(2n)$-th graded part).

Let $n$ be an arbitrary positive integer. Set
$$\wt{G}_n := \text{U}(2n, \F_{q^2}),\qquad \wt{\G}_n := \u(2n, \F_{q^2}).$$
For any $i, j\in\Z_{\ge 0}$, consider the maximal parabolic subgroup with Levi decomposition $\wt{P}_{i, j} = \wt{L}_{i, j}\,\wt{U}_{i, j}$ from \eqref{LU_unitary}. The Lie algebra of $\wt{L}_{i, j}$ is isomorphic to $\wt{\L}_{i, j}\cong\gl(i, \F_{q^2})\oplus\u(2j, \F_{q^2})$, and so $\CC(\wt{\L}_{i, j})\cong\CC(\gl(i, \F_{q^2}))\otimes\CC(\u(2j, \F_{q^2}))$. As a result, the parabolic restriction and induction can be interpreted as maps
\begin{gather*}
\Down^{\wt{\G}_n}_{\wt{\L}_{i, j}} : \CC(\u(2n, \F_{q^2})) \longrightarrow \CC(\gl(i, \F_{q^2}))\otimes\CC(\u(2j, \F_{q^2})),\\
\Up^{\wt{\G}_n}_{\wt{\L}_{i, j}} : \CC(\gl(i, \F_{q^2}))\otimes\CC(\u(2j, \F_{q^2})) \longrightarrow \CC(\u(2n, \F_{q^2})).
\end{gather*}
As in the previous section, these operations induce $\R$-linear graded maps
$$
\DeB: B_{q^2}\to A_{q^2}\otimes B_{q^2},\qquad \naB: A_{q^2}\otimes B_{q^2}\to B_{q^2}.
$$

The next theorem shows that the maps $\naB$, $\DeB$ equip $B_{q^2}$ with the structure of both an $A_{q^2}$-module and $A_{q^2}$-comodule, and the second theorem in this section explains the interaction between these maps.
The fact that we do not have a bialgebra and rather a bimodule boils down to the fact that the Levi subgroup $\wt{L}_{i, j}$ is isomorphic to a product of a general linear group $\GL(i, \F_{q^2})$ and an even unitary group $\UU(2j, \F_{q^2})$.

\begin{thm}\label{thm_twisted_0}
The map $\naB$ makes $B_{q^2}$ into a graded $A_{q^2}$-module and the map $\DeB$ makes it into a graded $A_{q^2}$-comodule. In other words:

(a) The maps $\naB$, $\DeB$ preserve grading.

(b) For any $a, a'\in A_{q^2}$, $b\in B_{q^2}$, we have
\begin{equation*}
\naB(\nabla\otimes\id)(a\otimes a'\otimes b) = \naB(\id\otimes\naB)(a\otimes a'\otimes b).
\end{equation*}

(c) For any $b\in B_{q^2}$, we have
\begin{equation*}
(\De\otimes\id)\DeB(b) = (\id\otimes\DeB)\DeB(b).
\end{equation*}
\end{thm}

\smallskip

It turns out that $B_{q^2}$ is not an $A_{q^2}$-bimodule, but actually certain \emph{twist} is needed.
For the next definition, it is helpful to use the following shorthand notations: for any $a_1, a_2, a_3\in A_{q^2}$, we denote
$$a_1 a_2 a_3 := \nabla(\nabla\otimes\id)(a_1\otimes a_2\otimes a_3) = \nabla(\id\otimes\nabla)(a_1\otimes a_2\otimes a_3).$$
Likewise, for any $a\in A_{q^2}$, $b\in B_{q^2}$, denote $a\cdot b := \wt{\nabla}(a\otimes b)$.

\begin{definition}\label{def_twist}
(i) For any $n\in\Z_{\geq 0}$, let $\om_n: \gl(n, \F_{q^2})\to \gl(n, \F_{q^2})$ be the involutive antilinear isomorphism $\om_n(X) := -J_nX^*J_n$, where $J_n$ is the $n\times n$ matrix \eqref{matrix_J}, i.e.
$$(\om_n(X))_{r, s} = -(X_{n+1-s,\ n+1-r})^q, \qquad X\in \gl(n, \F_{q^2}), \quad r,s\in\{1,\dots,n\}.$$
The maps $\omega_n$ induce linear maps $\CC(\gl(n, \F_{q^2}))\to\CC(\gl(n, \F_{q^2}))$, and therefore a graded linear map $\om: A_{q^2}\to A_{q^2}$.

(ii) Define the linear map $\odot : A_{q^2}\otimes (A_{q^2}\otimes B_{q^2}) \to (A_{q^2}\otimes B_{q^2})$ by
$$
\odot\left( a\otimes (a'\otimes b) \right) := \sum{(a_{(1)}\,\omega(a_{(2)})\, a')\otimes (a_{(3)}\cdot b)},\quad \forall\,a, a'\in A_{q^2},\ b\in B_{q^2},
$$
where
$$
(\Delta\otimes\id)\Delta(a) = (\id\otimes\Delta)\Delta(a) = \sum{a_{(1)}\otimes a_{(2)}\otimes a_{(3)}}.
$$
Moreover, denote
$$a\odot (a'\otimes b) := \odot\left( a\otimes (a'\otimes b) \right).$$
We call $\odot$ the \emph{twisted $A_{q^2}$-module structure on $A_{q^2}\otimes B_{q^2}$} (cf. van Leeuwen \cite{vL}).
\end{definition}

One can show that $\odot$ indeed turns $A_{q^2}\otimes B_{q^2}$ into an $A_{q^2}$-module.

%composition
%\begin{gather*}
%\odot := (\nabla\otimes\naB)\circ(\id_{A_{q^2}}\otimes P\otimes\id_{B_{q^2}})\circ(\xi\otimes \id_{A_{q^2}}\otimes \id_{A_{q^2}\otimes B_{q^2}})\circ(\De\otimes \id_{A_{q^2}\otimes B_{q^2}}),
%\end{gather*}
%where $P : A_{q^2}\otimes A_{q^2}\to A_{q^2}\otimes A_{q^2}$ is the linear map $P(a\otimes b) := b\otimes a$.

%The endomorphism $\xi: A_{q^2}\to A_{q^2}$ is the composition
%$$\xi := \nabla\circ (\id_{A_{q^2}}\otimes\,\om) \circ \De.$$

\begin{thm}[The twisted $A_{q^2}$-bimodule structure of $B_{q^2}$]\label{thm_twisted}
The map $\DeB$ intertwines the twisted action of $A_{q^2}$ on $A_{q^2}\otimes B_{q^2}$ with the left action of $A_{q^2}$ on $B_{q^2}$. In other words, for any $a\in A_{q^2}$, the following diagram is commutative:
\begin{center}
\begin{tikzcd}
B_{q^2}\arrow[r, "\DeB"] \arrow[d, swap, "b\mapsto \wt{\nabla}(a\otimes b)"]
& A_{q^2}\otimes B_{q^2} \arrow[d, "(a'\otimes b')\mapsto a\odot (a'\otimes b')"] \\
B_{q^2} \arrow[r, "\DeB"]
& A_{q^2}\otimes B_{q^2}
\end{tikzcd}
\end{center}
\end{thm}

\section{Positive harmonic functionals}\label{sec_8}

\subsection{On the bialgebra $A_q$}

Here, we relate positive harmonic functions on $\Gamma^{\GLB}$ and $\Gamma^{\GLB}_0 \sim \Y^{\HL}(q^{-1})$ to certain linear functionals on $A_q$.
Recall that $\T_n$ is the set of $\GL(n, \F_q)$-orbits of $\gl(n, \F_q)$, and $\T = \bigsqcup_{n\ge 0}{\T_n}$. For each $\bfla\in\T_n$, $n\in\Z_{\ge 0}$, denote the indicator function of the corresponding $\GL(n, \F_q)$-orbit by $\chi_{\bfla}\in\CC(\gl(n, \F_q))$.
Evidently, the set $\mathcal{X}_A := \{ \chi_{\bfla} \mid \bfla\in\T \}$ is a linear basis of $A_q$. Let $x_1\in\mathcal{X}_A$ be the indicator function of the ($\GL(1, \F_q)$-orbit of the) $1\times 1$ zero matrix $[0]\in\gl(1, \F_q)$.

\begin{definition}\label{FA}
Define $\calF(A_q)$ as the convex cone of linear functionals $\vp: A_q\to\R$ that satisfy:

$\bullet$ (Positivity) $\vp(x)\ge0,\ \forall\, x\in\mathcal{X}_A$, and

$\bullet$ (Harmonicity) $\vp(a) = \vp(\nabla(x_1\otimes a)),\ \forall\, a\in A_q$.

\smallskip

\noindent Moreover, let $\calF_0(A_q)\subset\calF(A_q)$ be the subcone of those functionals that satisfy the following vanishing condition:

$\bullet$ (Nilpotency) $\vp(x) = 0$, for all $x\in\mathcal{X}_A$ unless $x$ is a nilpotent orbit.
\end{definition}

\begin{proposition}\label{GL_functionals}
$\Harm_{\ge 0}(\Gamma^{\GLB})$ and $\calF(A_q)$ are isomorphic as convex cones.
Also, $\Harm_{\ge 0}(\Gamma^{\GLB}_0) \cong \Harm_{\ge 0}(\Y^{\HL}(q^{-1}))$ and $\calF_{\nil}(A_q)$ are isomorphic as convex cones.
\end{proposition}

\subsection{On the twisted $A_{q^2}$-bimodule $B_{q^2}$}

Recall that $\wt\T$ is the union of all $\UU(2n, \F_{q^2})$-orbits of $\u(2n, \F_{q^2})$, as $n$ ranges over $\Z_{\ge 0}$. For any $\wt{\bfla}\in\wt\T$, let $\chi_{\wt{\bfla}}$ be the indicator function of the unitary orbit $\wt\bfla$. The set $\mathcal{X}_B := \{ \chi_{\wt{\bfla}} \mid \wt{\bfla}\in\wt{\T} \}$ is a basis of $B_{q^2}$. Also, let $\wt{x}_1 : \CC(\gl(1, \F_{q^2}))\to\R$ be the indicator function of the $1\times 1$ zero matrix; observe that $\wt{x}_1\in A_{q^2}$ does not belong to $\mathcal{X}_B$.

\begin{definition}\label{FB}
$\wt{\calF}(B_{q^2})$ is the convex cone of linear functionals $\psi: B_{q^2}\to\R$ satisfying:

\smallskip

$\bullet$ (Positivity) $\psi(x)\ge0,\ \forall\, x\in\mathcal{X}_B$, and

$\bullet$ (Harmonicity) $\psi(b) = \psi\left( \wt{\nabla}(\wt{x}_1 \otimes b) \right),\ \forall \, b\in B_{q^2}$.

\noindent Furthermore, let $\wt{\calF}_0(B_{q^2})\subset\wt{\calF}(B_{q^2})$ be the subcone of those functionals which additionally satisfy:

$\bullet$ (Nilpotency) $\psi(x) = 0$, for all $x\in\mathcal{X}_B$ unless $x$ is a nilpotent unitary orbit.
\end{definition}

\begin{proposition}\label{U_functionals}
$\Harm_{\ge 0}(\Gamma^{\UB})$ and $\wt\calF(B_{q^2})$ are isomorphic as convex cones.
Also, $\Harm_{\ge 0}(\Gamma^{\UB}_0)\cong\Harm_{\ge 0}(\Y^{\HL}_{\E}(-q^{-1}))$ and $\wt\calF_{\nil}(B_{q^2})$ are isomorphic as convex cones.
\end{proposition}

\subsection{Reformulation of the main problem}

Recall that our initial motivation, in the language of branching graphs, was to answer Problem 2, i.e. to completely describe the convex cones $\Harm_{\ge 0}(\Gamma^{\GLB})$ and $\Harm_{\ge 0}(\Gamma^{\UB})$. As discussed before, a slight simplification is Problem 3, which asks for the characterization of the subcones $\Harm_{\ge 0}(\Gamma_0^{\GLB})$ and $\Harm_{\ge 0}(\Gamma_0^{\UB})$.
Based on Propositions \ref{GL_functionals} and \ref{U_functionals}, these problems can be reformulated as follows:

\begin{problem2a}
Study the convex cones $\calF(A_q),\, \wt\calF(B_{q^2})$, and characterize their extreme rays.
\end{problem2a}

\begin{problem3a}
Study the convex cones $\calF_{\nil}(A_q),\, \wt\calF_{\nil}(B_{q^2})$, and characterize their extreme rays.
\end{problem3a}

The structures of the convex cones $\calF(A_q),\, \calF_{\nil}(A_q)$ are completely understood. It remains to understand $\wt{\calF}(B_{q^2})$ and $\wt{\calF}_{\nil}(B_{q^2})$.
In the next section, we show how to use our knowledge on $\calF_{\nil}(A_{q^2})$ to produce a large infinite-parameter family of functionals in $\wt{\calF}_{\nil}(B_{q^2})$ from a single one.

\section{The mixing construction}\label{sec_9}

Motivated by Kerov, see \cite[Sec.~4]{GO} where his idea is elaborated, we introduce the following.

\begin{definition}[The mixing construction]
Let $\vp: A_{q^2}\to\R,\, \psi: B_{q^2}\to\R$ be any two real functionals, and let $s\in\R$ be any real number.
Then define a new linear functional
$$
\vp\star_s\psi: B_{q^2}\to\R
$$
by setting, for any homogeneous element $b'\in B_{q^2}$,
\begin{equation}\label{eqn_mixing}
(\vp\star_s\psi)(b'):=\sum_{a,b} s^{\deg a}(1-2s)^{\deg b}\,\vp(a)\psi(b),
\end{equation}
where $\wt{\De}(b')=\sum_{a,b} a\!\otimes\!b$ is any expansion such that $a\in A_{q^2}$ and $b\in B_{q^2}$ are homogeneous elements. In the sum \eqref{eqn_mixing}, we use the convention that $0^0 = 1$.
\end{definition}

\begin{thm}\label{kerov_construction}
Let\/ $0\le s\le \frac{1}{2}$ be arbitrary. Let $\varphi\in \calF(A_{q^2})$ and $\psi\in \wt{\calF}(B_{q^2})$ be positive harmonic functionals as in Definitions \ref{FA} and \ref{FB}, respectively. Then $\vp\star_s\psi$ belongs to $\wt{\calF}(B_{q^2})$. Moreover, if $\varphi\in\calF_0(A_{q^2})$ and $\psi\in\wt{\calF}_0(B_{q^2})$, then $\vp\star_s\psi$ belongs to $\wt{\calF}_0(B_{q^2})$.
\end{thm} 

\smallskip

We are now able to import our knowledge about $\calF_0(A_{q^2})$, in order to produce new families of functionals in $\wt\calF_0(B_{q^2})$.
Recall that the nilpotent orbits of $\T_n$ are parametrized by $\Y_n$; for any $\mu\in\Y_n$, let $\chi_\mu\in A_{q^2}$ be the indicator function of the corresponding orbit.
For example, $\chi_{(1)}$ is the indicator function of the ($GL(1, \F_{q^2})$-orbit of the) $1\times 1$ zero matrix. It is known from \cite{CO} that the functionals $\vp: A_{q^2}\to\R$ in the extreme rays of the cone $\calF_{\nil}(A_{q^2})$, normalized by $\vp(\chi_{(1)}) = 1$, are in bijection with
\begin{multline*}
\Omega(q^{-2}) := \{ \, \omega = (\alpha, \beta)\in(\R_{\ge 0})^{\infty}\times(\R_{\ge 0})^{\infty} \mid
\alpha_1\ge \alpha_2\ge \cdots\ge 0,\\
\beta_1\ge \beta_2\ge \cdots\ge 0,\quad \sum_{i\ge 1}{\alpha_i} + (1 - q^{-2})^{-1}\sum_{i\ge 1}{\beta_i} \le 1\}.
\end{multline*}
The infinite-dimensional space $\Omega(q^{-2})$ is a deformation of the \emph{Thoma simplex}, see \cite{BO}. If $\omega\in\Omega(q^{-2})$, denote the corresponding functional by $\vp_{\omega}\in\calF_0(A_{q^2})$; precise formulas for the values $\vp_\omega(\chi_\mu)$ can be derived from \cite[Thm. 4.18]{CO}.

Theorem \ref{kerov_construction} then shows that any single $\psi\in\wt\calF_{\nil}(B_{q^2})$ produces an infinite family $\{ \vp_{\omega}\star_s\psi \}_{\omega, s}$ of functionals in $\wt\calF_{\nil}(B_{q^2})$ parametrized by points $(\omega, s)\in\Omega(q^{-2})\times [0, \frac{1}{2}]$. As a first example, we consider the following functional $\psi_0: B_{q^2}\to\R$ (below $[0_{2n}]$ denotes the $(2n)\!\times\!(2n)$ zero matrix, which itself constitutes an orbit):
$$
\psi_0(\chi_{\wt{\bfla}}) :=
\begin{cases}
(q^2-1)^n\cdot\prod_{i=1}^{2n}{(q^i - (-1)^i)^{-1}}, & \text{if $\wt{\bfla}$ corresponds to $[0_{2n}]$ for some $n\ge0$};\\
0, & \text{otherwise}.
\end{cases}
$$

\begin{proposition}
$\psi_0$ belongs to $\wt{\calF}_{\nil}(B_{q^2})$.
\end{proposition}

In principle, different points $(\om, s)$ could give rise to the same functional $\vp_\om\star_s \psi_0$, but in fact this is impossible. Indeed, recall that nilpotent unitary orbits in $\wt\T_{2n}$ are parametrized by $\Y_{2n}$; denote the characteristic function corresponding to $\nu\in\Y_{2n}$ by $\wt{\chi}_\nu\in B_{q^2}$; then
\begin{equation*}
(\vp_\om \star_s \psi_0)( \wt{\chi}_{(1^{2n})} ) = \sum_{m=0}^n{ \frac{q^{3m^2 - 4mn}(q^2 - 1)^{n-m}}{\prod_{i=1}^{2n-2m}(q^i - (-1)^i)}\, s^m(1 - 2s)^{n-m}\ \vp_\omega(\chi_{(1^{m})})},\quad \forall\, n\in\Z_{\ge 0}.
\end{equation*}
These equalities (together with the fact that the quantities $\vp_\omega(\chi_{(1^{m})})$ uniquely determine the point $\omega\in\Omega(q^{-2})$) imply that the values of $\vp_\om\star_s \psi_0$ on the indicator functions of the singletons $[0_{2n}]$ uniquely determine the point $(\om, s)\in\Omega(q^{-2})\times [0, \frac{1}{2}]$.

\section*{Acknowledgements}

The research of the second author (G.~O.) was supported by the Russian Science Foundation, project 20-41-09009.

%\printbibliography


\begin{thebibliography}{9} 

\bibitem{B}
A. M. Borodin. The law of large numbers and the central limit theorem for the Jordan normal form of large triangular matrices over a finite field. Journal of Mathematical Sciences (New York) 96, no. 5 (1999), 3455--3471.

\bibitem{BO}
A. Borodin and G. Olshanski. Representations of the infinite symmetric group. Cambridge Studies in Advanced Mathematics 160. Cambridge University Press, 2017.

\bibitem{CO}
C. Cuenca and G. Olshanski. Infinite-dimensional groups over finite fields and Hall-Littlewood symmetric functions. Advances in Mathematics, vol. 395 (2022), 108087.

\bibitem{CR}
C. W. Curtis and I. Reiner. Representation theory of finite groups and associative algebras. Amer. Math. Soc. Chelsea Publishing. 1962. 

\bibitem{DM}
F. Digne and J. Michel. Representations of finite groups of Lie type, 2nd ed. London Mathematical Society Student Texts 95.  Cambridge University Press, 2020.

\bibitem{F}
J. Fulman, Random matrix theory over finite fields, Bull. Amer. Math. Soc. 39 (2002) 51--85.

\bibitem{GO}
A. Gnedin and G. Olshanski. Coherent permutations with descent statistic and the boundary problem for the graph of zigzag diagrams. International Mathematics Research Notices 2006 (2006), paper 51968.

\bibitem{GKV}
V. Gorin, S. Kerov and A. Vershik. Finite traces and representations of the group of infinite matrices over a finite field. Advances in Math. 254 (2014), 331--395.

\bibitem{GOl}
V. Gorin and G. Olshanski. A quantization of the harmonic analysis on the infinite-dimensional unitary group. Journal of Functional Analysis 270, no. 1 (2016), 375--418.

\bibitem{GR}
D. Grinberg and V. Reiner. Hopf Algebras in Combinatorics (2014). Preprint, arXiv:1409.8356.

\bibitem{Ke0}
S. V. Kerov. Combinatorial examples in the theory of AF-algebras.  J. Math. Sci. (New York)  59:5 (1992), 1063--1071.

\bibitem{Ke}
S. V. Kerov. Generalized Hall-Littlewood symmetric functions and orthogonal polynomials. Representation theory and dynamical systems, Adv. Soviet Math., 9, Amer. Math. Soc., Providence, RI (1992), 67--94.

\bibitem{KOO}
S. Kerov, A. Okounkov and G. Olshanski. The boundary of the Young graph with Jack edge multiplicities. International Mathematics Research Notices 1998, Vol. 4 (1998), 173--199.

\bibitem{K}
A. A. Kirillov. Variations on the triangular theme. In: Lie Groups and Lie Algebras: E. B. Dynkin's Seminar. (Amer. Math. Soc. Transl. Ser. 2, vol. 169). AMS,1995, 43--73.

\bibitem{vL}
M. A. A. van Leeuwen. An application of Hopf-algebra techniques to representations of finite classical groups. Journal of Algebra 140, issue 1 (1991), 210--246.

\bibitem{Le}
G. Lehrer. The space of invariant functions on a finite Lie algebra. Transactions of the American Mathematical Society 348, no. 1 (1996), 31--50.

\bibitem{Mac}
I. G. Macdonald. Symmetric functions and Hall polynomials. Oxford university press, 1998.

\bibitem{Mackey}
 G. W. Mackey. On induced representations of groups. Amer. J. Math. 73 (1951), 576--592.

\bibitem{Mat}
K. Matveev. Macdonald-positive specializations of the algebra of symmetric functions: Proof of the Kerov conjecture. Annals of Mathematics 189, Issue 1 (2019), 277--316.

\bibitem{OO1}
A. Okounkov and G. Olshanski. Asymptotics of Jack polynomials as the number of variables goes to infinity. International Mathematics Research Notices 1998, no. 13 (1998), pp. 641--682.

\bibitem{OO2}
A. Okounkov and G. Olshanski. Limits of BC-type orthogonal polynomials as the number of variables goes to infinity. In: Jack, Hall-Littlewood and Macdonald Polynomials, American Mathematical Society Contemporary Mathematics Series 417 (2006), pp. 281--318.

\bibitem{OV}
G. Olshanski and A. Vershik. Ergodic unitarily invariant measures on the space of infinite Hermitian matrices. Translations of the American Mathematical Society-Series 2 175 (1996), 137--176.

\bibitem{Sc}
O. Schiffmann. Lectures on Hall algebras. In: Geometric methods in representation theory. II, volume 24 of S\'emin. Congr., pages 1--141. Soc. Math. France, Paris, 2012.

\bibitem{Sp}
T. A. Springer. The Steinberg function of a finite Lie algebra. Inventiones mathematicae 58, no. 3 (1980), 211--215.

\bibitem{T}
E. Thoma. Die unzerlegbaren, positive-definiten Klassenfunktionen der
abz\"ahlbar unendlichen, symmetrischen Gruppe. Math. Zeitschr. 85
(1964), 40--61.

\bibitem{VK1}
A. M. Vershik and S. V. Kerov. Asymptotic theory of characters of the
symmetric group. Funct. Anal. Appl. 15 (1981), no. 4, 246--255.

\bibitem{VK2}
A. M. Vershik and S. V. Kerov. Characters and factor representations of
the infinite unitary group. Doklady AN SSSR 267 (1982), no. 2, 272--276
(Russian); English translation: Soviet Math. Doklady 26 (1982), 570--574.

\bibitem{V82}
A. M. Vershik. Symmetric Functions and K–functor, in the appendix of translation editor to G. James, Representation theory of symmetric groups (Russian edition). Moscow, Mir, 1982.

\bibitem{VK98}
A. M. Vershik and S. V. Kerov. On a infinite-dimensional group over a finite field. Functional Analysis and Its Applications, Vol. 32, no. 3 (1998), 147--152.

\bibitem{VK07}
A. M. Vershik and S. V. Kerov. Four drafts of the representation theory of the group of infinite matrices over a finite field. Journal of Mathematical Sciences, Vol. 147, no. 6 (2007), 7129--7144.

\bibitem{V}
D. Voiculescu, Repr\'esentations factorielles de type {\rm II}${}_1$ de
$U(\infty)$. J. Math. Pures et Appl. 55 (1976), 1--20.

\bibitem{Z}
A. V. Zelevinsky. Representations of finite classical groups: a Hopf algebra approach. Vol. 869. Springer, 2006.
\end{thebibliography}
\end{document}